# ADAPTIVE HAUSDORFF ESTIMATION OF DENSITY LEVEL SETS

By Aarti Singh,[1] Clayton Scott and Robert Nowak[1]

*University of Wisconsin–Madison, University of Michigan–Ann Arbor and University of Wisconsin–Madison*

Consider the problem of estimating the $\gamma$-level set $G^*_\gamma = \{x : f(x) \geq \gamma\}$ of an unknown $d$-dimensional density function $f$ based on $n$ independent observations $X_1, \ldots, X_n$ from the density. This problem has been addressed under global error criteria related to the symmetric set difference. However, in certain applications a spatially uniform mode of convergence is desirable to ensure that the estimated set is close to the target set everywhere. The Hausdorff error criterion provides this degree of uniformity and, hence, is more appropriate in such situations. It is known that the minimax optimal rate of error convergence for the Hausdorff metric is $(n/\log n)^{-1/(d+2\alpha)}$ for level sets with boundaries that have a Lipschitz functional form, where the parameter $\alpha$ characterizes the regularity of the density around the level of interest. However, the estimators proposed in previous work are nonadaptive to the density regularity and require knowledge of the parameter $\alpha$. Furthermore, previously developed estimators achieve the minimax optimal rate for rather restricted classes of sets (e.g., the boundary fragment and star-shaped sets) that effectively reduce the set estimation problem to a function estimation problem. This characterization precludes level sets with multiple connected components, which are fundamental to many applications. This paper presents a fully data-driven procedure that is adaptive to unknown regularity conditions and achieves near minimax optimal Hausdorff error control for a class of density level sets with very general shapes and multiple connected components.

**1. Introduction.** Level sets provide useful summaries of a function for many applications including clustering [6, 8, 21], anomaly detection [16, 20, 24], functional neuroimaging [12, 25], bioinformatics [27], digital elevation

Received November 2007; revised August 2008.

[1]Supported in part by NSF Grants ECS-05-29381, CCF-03-53079 and CCR-03-50213.

mapping [19, 26] and environmental monitoring [22]. In practice, however, the function itself is unknown a priori, and only a finite number of observations related to $f$ are available. In this paper, we focus on the density level set problem; extensions to general regression level set estimation should be possible using a similar approach, but they are beyond the scope of this paper. Let $X_1, \ldots, X_n$ be independent, identically distributed observations drawn from an unknown probability measure $P$, having density $f$ with respect to the Lebesgue measure and defined on the domain $\mathcal{X} \subseteq \mathbb{R}^d$. Given a desired density level $\gamma$, consider the $\gamma$-level set of the density $f$

$$G^*_\gamma := \{x \in \mathcal{X} : f(x) \geq \gamma\}.$$

The goal of the density level set estimation problem is to generate an estimate $\widehat{G}$ of the level set based on the $n$ observations $\{X_i\}_{i=1}^n$, such that the error between the estimator $\widehat{G}$ and the target set $G^*_\gamma$, as assessed by some performance measure which gauges the closeness of the two sets, is small.

Most literature available on level set estimation methods [9, 13, 14, 15, 16, 20, 23, 26] considers error measures related to the symmetric set difference, $G_1 \Delta G_2 = (G_1 \setminus G_2) \cup (G_2 \setminus G_1)$. However, level set methods based on a measure of the symmetric difference error may produce estimates that veer greatly from the desired level set at certain places, since the symmetric difference is a global measure of *average* closeness between two sets. Some applications may need a more local or spatially uniform error measure as provided by the Hausdorff metric, for example, to preserve topological properties of the level set as in clustering [6, 8, 21] or ensure robustness to outliers in level set-based anomaly detection [16, 20, 24] and data ranking [11]. The Hausdorff error metric is defined as follows between two nonempty sets:

$$d_\infty(G_1, G_2) = \max\Big\{\sup_{x \in G_2} \rho(x, G_1), \sup_{x \in G_1} \rho(x, G_2)\Big\},$$

where $\rho(x, G) = \inf_{y \in G} \|x - y\|$, the smallest Euclidean distance of a point in $G$ to the point $x$. If $G_1$ or $G_2$ is empty, then let $d_\infty(G_1, G_2)$ be defined as the largest distance between any two points in the domain. Control of this error measure provides a uniform mode of convergence, as it implies control of the deviation of a single point from the desired set. A symmetric set difference-based estimator may not provide such a uniform control as it is easy to see that a set estimate can have a very small measure of symmetric difference error but large Hausdorff error. Conversely, as long as the set boundary is not space filling and the domain is bounded, small Hausdorff error implies small symmetric-difference measure.

Existing results pertaining to nonparametric level set estimation using the Hausdorff metric [2, 9, 23] focus on rather restrictive classes of level sets (e.g., the boundary fragment and star-shaped set classes). These restrictions, which effectively reduce the set estimation problem to a boundary function estimation problem (in rectangular or polar coordinates, resp.), are typically not met in practical applications. In particular, the characterization of level set estimation as a boundary function estimation problem requires prior knowledge of a reference coordinate or interior point (in rectangular or polar coordinates, resp.) and precludes level sets with multiple connected components. Moreover, the estimation techniques proposed in [2, 9, 23] require precise knowledge of the regularity of the density (quantified by the parameter $\alpha$, to be defined below) in the vicinity of the desired level in order to achieve minimax optimal rates of convergence. Such prior knowledge is unavailable in most practical applications. Recently, a plug-in method based on sup-norm density estimation was put forth in [3] that can handle more general classes than boundary fragments or star-shaped sets. However, sup-norm density estimation requires the density to satisfy global smoothness assumptions. Also, the method only deals with a special case of the density regularity condition considered in this paper ($\alpha = 1$) and is therefore not adaptive to unknown density regularity.

In this paper, we propose a plug-in procedure based on a regular histogram partition that can adaptively achieve minimax optimal rates of Hausdorff error convergence over a broad class of level sets with very general shapes and multiple connected components, without assuming a priori knowledge of the density regularity parameter $\alpha$. Adaptivity is achieved by a new data-driven procedure for selecting the histogram resolution. The procedure bears some similarity to Lepski-type methods [10], as further discussed in Section 3.2. However, our procedure is specifically designed for the level set estimation problem and only requires *local* regularity of the density in the vicinity of the desired level. A shorter version of this paper appeared in [17]; however, it relies on more stringent assumptions on the class of level sets under consideration. In this paper, we generalize the class of level sets to allow for spatial variations in the density regularity along the level set boundary, and we also discuss extensions to support set estimation and discontinuity in the density at all points around the level of interest.

The paper is organized as follows. Section 2 states our basic assumptions which allow Hausdorff accurate level set estimation and presents a minimax lower bound on the Hausdorff performance of any level set estimator for the class of densities under consideration. In Section 3, we present the proposed histogram-based approach to Hausdorff accurate level set estimation. In Section 3.1, we show that the proposed estimator can achieve the minimax optimal rate of convergence given knowledge of the density regularity parameter $\alpha$, and Section 3.2 extends the estimator to achieve adaptivity to unknown density regularity. We also comment on extensions that address discontinuity in the density at the level of interest and support set estimation. Concluding remarks are given in Section 4 and the Appendices contain proofs of the main results.

**2. Density assumptions.** We assume that the domain of the density $f$ is the unit hypercube in $d$ dimensions, that is, $\mathcal{X} = [0,1]^d$. Extensions to other compact domains are straightforward. Furthermore, the density is assumed to be bounded with range $[0, f_{\max}]$, though we do not assume knowledge of $f_{\max}$. Controlling the Hausdorff accuracy of level set estimates requires some smoothness assumptions on the density and the level set boundary, which are stated below. Before that, we introduce the following definitions:

- *$\varepsilon$-ball*: An $\varepsilon$-ball centered at a point $x \in \mathcal{X}$ is defined as

$$B(x, \varepsilon) = \{y \in \mathcal{X} : \|x - y\| \leq \varepsilon\}.$$

  Here $\|\cdot\|$ denotes the Euclidean distance.
- *Inner $\varepsilon$-cover*: An inner $\varepsilon$-cover of a set $G \subseteq \mathcal{X}$ is defined as the union of all $\varepsilon$-balls contained in $G$. Formally,

$$\mathcal{I}_\varepsilon(G) = \bigcup_{x \,:\, B(x,\varepsilon) \subseteq G} B(x, \varepsilon).$$

We are now ready to state the assumptions. The first one characterizes the relationship between distances and changes in density, and the second one is a topological assumption on the level set boundary that essentially generalizes the notion of Lipschitz functions to closed hypersurfaces.

[A] *Local density regularity.* The density is $\alpha$-regular around the $\gamma$-level set, $0 < \alpha < \infty$ and $0 < \gamma < f_{\max}$, if:
  [A1] there exist constants $C_1, \delta_1 > 0$ such that for all $x \in \mathcal{X}$ with $|f(x) - \gamma| \leq \delta_1$,

$$|f(x) - \gamma| \geq C_1 \rho(x, \partial G_\gamma^*)^\alpha,$$

  where $\partial G_\gamma^*$ denotes the boundary of the true level set $G_\gamma^*$.
  [A2] there exist constants $C_2, \delta_2 > 0$ and $x_0 \in \partial G_\gamma^*$ such that for all $x \in B(x_0, \delta_2)$,

$$|f(x) - \gamma| \leq C_2 \rho(x, \partial G_\gamma^*)^\alpha.$$

This condition characterizes the behavior of the density around the level $\gamma$. Assumption [A1] states that the density cannot be arbitrarily "flat" around the level, and changes as at least the $\alpha$th power of the distance from the level set boundary. Assumption [A2] states that there exists a fixed neighborhood around some point on the boundary where the density changes no faster than the $\alpha$th power of the distance from the level set boundary. The latter condition is only required for adaptivity, as we discuss later. The regularity parameter $\alpha$ determines the rate of error convergence for level set estimation. Accurate estimation is more difficult at levels where the density is relatively flat (large $\alpha$), as intuition would suggest. It is important to point out that in this paper we do not assume knowledge of $\alpha$, unlike previous investigations into Hausdorff accurate level set estimation [2, 3, 9, 23]. Therefore, here the assumption simply states that there is a relationship between distance and density level, but the precise nature of the relationship is unknown. In Section 3, we briefly discuss extensions to address the case $\alpha = 0$ which corresponds to discontinuity in the density at all points around the level set boundary and the case $\gamma = 0$ which corresponds to support set estimation.

[B] *Level set regularity.* There exist constants $\varepsilon_o > 0$ and $C_3 > 0$ such that for all $\varepsilon \leq \varepsilon_o$, $\mathcal{I}_\varepsilon(G^*_\gamma) \neq \varnothing$, and for all $x \in \partial G^*_\gamma$, $\rho(x, \mathcal{I}_\varepsilon(G^*_\gamma)) \leq C_3 \varepsilon$. This assumption implies that the level set is not arbitrarily narrow anywhere. It precludes space-filling boundaries and features like cusps, arbitrarily thin ribbons and isolated connected components of arbitrarily small size. This condition is necessary since arbitrarily small features cannot be detected and resolved from a finite sample.

For a fixed set of positive numbers $C_1$, $C_2$, $C_3$, $\varepsilon_0$, $\delta_1$, $\delta_2$, $f_{\max}$, $\gamma < f_{\max}$, $d$ and $\alpha$, we consider the following classes of densities.

DEFINITION 1. $\mathcal{F}^*_1(\alpha)$ denotes the class of densities satisfying assumptions [A1] and [B].

DEFINITION 2. $\mathcal{F}^*_2(\alpha)$ denotes the class of densities satisfying assumptions [A1], [A2] and [B].

The dependence on other parameters is omitted as these do not influence the minimax optimal rate of convergence (except for the dimension $d$). In the paper, we present a method that provides minimax optimal rates of convergence for the class $\mathcal{F}^*_1(\alpha)$, given knowledge of the density regularity parameter $\alpha$. We also extend the method to achieve adaptivity to $\alpha$ for the class $\mathcal{F}^*_2(\alpha)$, while preserving the minimax optimal performance.

Assumption [A] is similar to the one employed in [2, 23], except that the upper bound assumption on the density deviation in [2, 23] holds provided that the set $\{x : |f(x) - \gamma| \leq \delta_1\}$ is nonempty. This implies that the densities either jump across the level $\gamma$ at any point on the level set boundary (i.e., the deviation is greater than $\delta_1$) or change exactly as the $\alpha$th power of the distance from the boundary. Our formulation allows for densities with regularities that vary spatially along the level set boundary—it requires that the density changes no slower than the $\alpha$th power of the distance from the boundary, except in a fixed neighborhood of one point where the density changes exactly as the $\alpha$th power of the distance from the boundary. While the formulation in [2, 23] requires the upper bound on the density deviation to hold for at least one point on the boundary, our assumption [A2] requires the upper bound to hold for a fixed neighborhood around at least one point on the boundary. This is necessary for adaptivity since a procedure cannot sense the regularity as characterized by $\alpha$ if the regularity only holds in an arbitrarily small region. Assumption [B] implies that the boundary looks locally like a Lipschitz function and allows for level sets with multiple connected components and arbitrary locations. Thus, these restrictions are quite mild and less restrictive than those considered in the previous literature on Hausdorff accurate level set estimation. In fact, assumption [B] is satisfied by a Lipschitz boundary fragment or star-shaped set as considered in [2, 9, 23], as the following lemma states; please refer to [18] for a formal proof.

LEMMA 1. *Consider the $\gamma$ level set $G^*_\gamma$ of a density $f \in \mathcal{F}_{SL}(\alpha)$, where $\mathcal{F}_{SL}(\alpha)$ denotes the class of $\alpha$-regular densities with Lipschitz star-shaped level sets as defined in [23]. Then, $G^*_\gamma$ satisfies the level set regularity assumption* [B].

In Theorem 4 of [23], Tsybakov establishes a minimax lower bound of $(n/\log n)^{-1/(d+2\alpha)}$ for the class of Lipschitz star-shaped sets, which, per Lemma 1, also satisfy assumption [B]. His proof uses Fano's lemma to derive the lower bound for a discrete subset of densities from this class. It is easy to see that the discrete subset of densities used in his construction also satisfy our form of assumption [A]. Hence, the same minimax lower bound holds for the classes $\mathcal{F}^*_1(\alpha)$ and $\mathcal{F}^*_2(\alpha)$ under consideration as well, and we have the following proposition. Here $\mathbb{E}$ denotes expectation with respect to the random data sample.

PROPOSITION 1. *There exists $c > 0$ such that, for large enough $n$,*

$$\inf_{G_n} \sup_{f \in \mathcal{F}^*_1(\alpha)} \mathbb{E}[d_\infty(G_n, G^*_\gamma)] \geq \inf_{G_n} \sup_{f \in \mathcal{F}^*_2(\alpha)} \mathbb{E}[d_\infty(G_n, G^*_\gamma)] \geq c\left(\frac{n}{\log n}\right)^{-1/(d+2\alpha)}.$$

*The* inf *is taken over all set estimators $G_n$ based on the $n$ observations.*

**3. Hausdorff accurate level set estimation using histograms.** Direct Hausdorff estimation is challenging as there exists no natural empirical measure that can be used to gauge the Hausdorff error of an estimate. However, the density regularity assumption [A] suggests that Hausdorff control over the level set estimate can be obtained indirectly by controlling the density deviation error rather than the distance deviation. Thus, we propose a plug-in level set estimator that is based on an empirical density estimator, the regular histogram.

Let $\mathcal{A}_j$ denote the collection of cells in a regular partition of $[0,1]^d$ into hypercubes of dyadic sidelength $2^{-j}$, where $j$ is a nonnegative integer. The level set estimate at this resolution is given as

$$\widehat{G}_j = \bigcup_{A \in \mathcal{A}_j \,:\, \widehat{f}(A) \geq \gamma} A. \tag{1}$$

Here $\widehat{f}(A) = \widehat{P}(A)/\mu(A)$, where $\widehat{P}(A) = \frac{1}{n} \sum_{i=1}^n \mathbf{1}_{\{X_i \in A\}}$ denotes the empirical probability of an observation occurring in $A$, and $\mu$ is the Lebesgue measure.

3.1. *A priori knowledge of local density regularity.* The appropriate resolution for accurate level set estimation depends on the density regularity, as characterized by $\alpha$, near the level of interest. If the density varies sharply near the level of interest (small $\alpha$), then accurate estimation is easier and a fine resolution suffices. Identifying the level set is more difficult if the density is very flat (large $\alpha$) and, hence, a lower resolution (more averaging) is required. Our first result shows that if the local density regularity parameter $\alpha$ is known, then the correct resolution for Hausdorff accurate level set estimation can be chosen (as in [2, 23]), and the corresponding estimator of (1) achieves near minimax optimal rate over the class of densities given by $\mathcal{F}_1^*(\alpha)$. Notice that even though the proposed method is a plug-in level set estimator based on a histogram density estimate, the histogram resolution is chosen to specifically target the level set problem and is not optimized for density estimation. Thus, we do not require that the density exhibits some smoothness at all points in the domain. We introduce the notation $a_n \asymp b_n$ to denote $a_n = O(b_n)$ and $b_n = O(a_n)$.

THEOREM 1. *Assume that the local density regularity $\alpha$ is known. Pick resolution $j \equiv j(n)$ such that $2^{-j} \asymp s_n (n/\log n)^{-1/(d+2\alpha)}$, where $s_n$ is a monotone diverging sequence. Then,*

$$\sup_{f \in \mathcal{F}_1^*(\alpha)} \mathbb{E}[d_\infty(\widehat{G}_j, G_\gamma^*)] \leq C s_n \left(\frac{n}{\log n}\right)^{-1/(d+2\alpha)}$$

*for all $n$, where $C \equiv C(C_1, C_3, \varepsilon_o, f_{\max}, \delta_1, d, \alpha) > 0$ is a constant.*

The proof is given in Appendix A and relies on two key facts. First, the density regularity assumption [A1] implies that the distance of any point in the level set estimate is controlled by its deviation from the level of interest $\gamma$. Therefore, with high probability, only the cells near the boundary are erroneously included or excluded in the level set estimate. Second, the level set boundary does not have very narrow features—features that cannot be detected by a finite sample—and is locally Lipschitz as per assumption [B].

This implies that the erroneous cells are not too far from the nonerroneous cells. Using these arguments, it is shown that the Hausdorff error scales as the histogram cell sidelength.

Theorem 1 provides an upper bound on the Hausdorff error of our estimate. If $s_n$ is slowly diverging, for example, if $s_n = (\log n)^\varepsilon$ where $\varepsilon > 0$, this upper bound agrees with the minimax lower bound of Proposition 1 up to a $(\log n)^\varepsilon$ factor. Hence, the proposed estimator can achieve near minimax optimal rates, given knowledge of the density regularity. We would like to point out that if the parameter $\delta_1$ characterizing assumption [A] and the density bound $f_{\max}$ are also known, then the appropriate resolution can be chosen as $j = \lfloor \log_2(c^{-1}(n/\log n)^{1/(d+2\alpha)}) \rfloor$, where the constant $c \equiv c(\delta_1, f_{\max})$. With this choice, the optimal sidelength scales as $2^{-j} \asymp (n/\log n)^{-1/(d+2\alpha)}$, and the estimator $\widehat{G}_j$ exactly achieves the minimax optimal rate.

REMARK 1. A dyadic sidelength is not necessary for Theorem 1 to hold, however the adaptive procedure described below is based on a search over dyadic resolutions. Thus, to present a unified analysis, we consider a dyadic sidelength here as well.

3.2. *Adapting to unknown local density regularity.* In this section, we present a procedure that automatically selects the appropriate resolution in a purely data-driven way without assuming prior knowledge of $\alpha$. The proposed procedure is a complexity regularization approach that is reminiscent of Lepski-type methods for function estimation [10], which are spatially adaptive bandwidth selectors. In Lepski methods, the appropriate bandwidth at a point is determined as the largest bandwidth for which the estimate does not deviate significantly from estimates generated at finer resolutions. Our procedure is similar in spirit, however it is tailored specifically for the level set problem; hence, the chosen resolution at any point depends only on the local regularity of the density around the level of interest.

The histogram resolution search is focused on regular partitions of dyadic sidelength $2^{-j}$, $j \in \{0, 1, \ldots, J\}$. The choice of $J$ will be specified below. Since the selected resolution needs to be adapted to the local regularity of the density around the level of interest, we introduce the following vernier:

$$\mathcal{V}_{\gamma,j} = \min_{A \in \mathcal{A}_j} \max_{A' \in \mathcal{A}_{j'} \cap A} |\gamma - \bar{f}(A')|.$$

Here $\bar{f}(A) = P(A)/\mu(A)$, $j' = \lfloor j + \log_2 s_n \rfloor$, where $s_n$ is a slowly diverging monotone sequence, for example, $\log n$, $\log \log n$, etc., and $\mathcal{A}_{j'} \cap A$ denotes the collection of subcells with sidelength $2^{-j'} \in [2^{-j}/s_n, 2^{-j+1}/s_n)$ within the cell $A$. Observe that the vernier value is determined by a cell $A \in \mathcal{A}_j$ that intersects the boundary $\partial G^*_\gamma$. By evaluating the deviation in average density from level $\gamma$ within subcells of $A$, the vernier indicates whether or not the density in cell $A$ is uniformly close to $\gamma$. Thus, the vernier is sensitive to the local density regularity in the vicinity of the desired level and leads to selection of the appropriate resolution adapted to the unknown density regularity parameter $\alpha$, as we will show in Theorem 2.

Since $\mathcal{V}_{\gamma,j}$ requires knowledge of the unknown probability measure, we must work with the empirical version, defined analogously as

$$\widehat{\mathcal{V}}_{\gamma,j} = \min_{A \in \mathcal{A}_j} \max_{A' \in \mathcal{A}_{j'} \cap A} |\gamma - \widehat{f}(A')|.$$

The empirical vernier $\widehat{\mathcal{V}}_{\gamma,j}$ is balanced by a penalty term

$$\Psi_{j'} := \max_{A \in \mathcal{A}_{j'}} \sqrt{8 \frac{\log(2^{j'(d+1)} 16/\delta)}{n\mu(A)} \max\left(\widehat{f}(A), 8 \frac{\log(2^{j'(d+1)} 16/\delta)}{n\mu(A)}\right)},$$

where $0 < \delta < 1$ is a confidence parameter, and $\mu(A) = 2^{-j'd}$. Notice that the penalty is computable from the given observations. The precise form of $\Psi$ is chosen to bound the deviation between true and empirical vernier with high probability (refer to Corollary B.1 for a formal proof). The final level set estimate is given by

$$\widehat{G} = \widehat{G}_{\widehat{j}}, \tag{2}$$

where

$$\widehat{j} = \arg \min_{0 \leq j \leq J} \{\widehat{\mathcal{V}}_{\gamma,j} + \Psi_{j'}\}. \tag{3}$$

Observe that the value of the vernier decreases with increasing resolution as better approximations to the true level are available. On the other hand, the penalty is designed to increase with resolution to penalize high complexity estimates that might overfit the given sample of data. Thus, the above procedure chooses the appropriate resolution automatically by balancing these two terms. The following theorem characterizes the performance of the proposed complexity penalized procedure.

THEOREM 2. *Pick $J \equiv J(n)$ such that $2^{-J} \asymp s_n(n/\log n)^{-1/d}$, where $s_n$ is a monotone diverging sequence. Let $\widehat{j}$ denote the resolution chosen by the complexity penalized method as given by (3) and $\widehat{G}$ denote the final estimate of (2). Then, with probability at least $1 - 2/n$, for all densities in the class $\mathcal{F}_2^*(\alpha)$,*

$$c_1 s_n^{d/(d+2\alpha)} \left(\frac{n}{\log n}\right)^{-1/(d+2\alpha)} \leq 2^{-\widehat{j}} \leq c_2 s_n s_n^{d/(d+2\alpha)} \left(\frac{n}{\log n}\right)^{-1/(d+2\alpha)}$$

*for $n$ large enough [so that $s_n > c(C_3, \varepsilon_o, d)$], where $c_1, c_2 > 0$ are constants. In addition,*

$$\sup_{f \in \mathcal{F}_2^*(\alpha)} \mathbb{E}[d_\infty(\widehat{G}, G_\gamma^*)] \leq C s_n^2 \left(\frac{n}{\log n}\right)^{-1/(d+2\alpha)}$$

*for all* $n$, *where* $C \equiv C(C_1, C_2, C_3, \varepsilon_o, f_{\max}, \delta_1, \delta_2, d, \alpha) > 0$ *is a constant.*

The proof is given in Appendix B. Observe that the maximum resolution $2^J \asymp s_n^{-1}(n/\log n)^{1/d}$ depends only on $n$ and allows the optimal resolution for any $\alpha$ to lie in the search space. By appropriate choice of $s_n$, for example, $s_n = (\log n)^{\varepsilon/2}$ with $\varepsilon$ a small number $> 0$, the bound of Theorem 2 matches the minimax lower bound of Proposition 1, except for an additional $(\log n)^\varepsilon$ factor. Hence, our method *adaptively* achieves near minimax optimal rates of convergence for the class $\mathcal{F}_2^*(\alpha)$.

REMARK 2. The case $\alpha = 0$ corresponds to jump in the density across the level $\gamma$, at all points along the level set boundary. The adaptive estimator can be extended to handle the complete range $0 \leq \alpha < \infty$ by a slight modification of the vernier

$$\mathcal{V}_{\gamma,j} = 2^{-j'/2} \min_{A \in \mathcal{A}_j} \max_{A' \in \mathcal{A}_{j'} \cap A} |\gamma - \bar{f}(A')|.$$

This makes the vernier sensitive to the resolution even for the jump case and biases a vernier minimizer toward finer resolutions. The exact form of the modification arises from technical considerations and is somewhat nonintuitive. Hence, we omitted the jump case in our earlier analysis to keep the presentation simple. The penalty also needs to be scaled by a factor of $2^{-j'/2}$, to ensures that balancing the vernier and penalty leads to the appropriate resolution for the whole range of the regularity parameter $0 \leq \alpha < \infty$. Please refer to [18] for a detailed proof.

REMARK 3. Under a measure of the symmetric difference error, it is known that support set estimation, that is, learning the set $G_0^* := \{x : f(x) > 0\}$, is easier than level set estimation, except for the case $\alpha = 0$ (see [7, 23]). The same holds for Hausdorff error and the minimax rate of convergence can be shown to be $(n/\log n)^{-1/(d+\alpha)}$ [18]. The minimax lower bound follows along the lines of the minimax lower bound in [23] for level set estimation $(\gamma > 0)$. This rate can be achieved by the following plug-in histogram estimator:

$$\widehat{G}_{0,j} = \bigcup_{A \in \mathcal{A}_j \,:\, \widehat{f}(A) > 0} A.$$

The analysis requires a modified theoretical analysis using Bernstein inequalities rather than the relative VC inequalities we use in the proofs of Theorems 1 and 2 for level set estimation. Formal proofs for support set estimation are given in [18].

**4. Conclusions.** In this paper, we developed a Hausdorff accurate level set estimation method that is adaptive to unknown density regularity and achieves nearly minimax optimal rates of error convergence over a more general class of level sets than considered in previous literature. The vernier provides the key to achieve adaptivity while requiring only local regularity of the density in the vicinity of the desired level. We also discussed extensions of the proposed estimator to address discontinuity in the density around the level of interest and support set estimation.

While this paper considers level sets with locally Lipschitz boundaries, extensions to additional boundary smoothness (e.g., Hölder regularity $>$ 1) may be possible in the proposed framework using techniques such as wedgelets [5] or curvelets [1]. The earlier work on Hausdorff accurate level set estimation [2, 9, 23] does address higher smoothness of the boundary, but that follows as a straightforward consequence of assuming a functional form for the boundary. Also, we have only addressed the density level set problem in this paper. Extensions to general regression level set estimation should be possible using a similar approach.

The results of this paper indicate that a regular, spatially nonadaptive partition suffices for minimax optimal Hausdorff accurate level set estimation. However, in practice, a spatially adapted partition can provide better performance than a uniform partition. This is because nonuniform partitions can adapt to the spatial variations in density regularity to yield better estimate of the boundary where the density changes sharply, even though the Hausdorff error is dominated by the accuracy in regions where the density is relatively flat at the level of interest. Thus, it is of interest to develop spatially adapted estimators. This might be possible by developing a tree-based approach or a modified Lepski method, and it is the subject of current research.

## APPENDIX A: PROOF OF THEOREM 1

Before proceeding to the proof of Theorem 1, we establish three lemmas that will be used both in this proof and in the proof of Theorem 2. The first lemma bounds the deviation of true and empirical density averages. The choice of penalty used to achieve adaptivity is motivated by this relation.

LEMMA A.1. *Consider $0 < \delta < 1$. With probability at least $1 - \delta$, the following is true for all $j \geq 0$*

$$\max_{A \in \mathcal{A}_j} |\bar{f}(A) - \widehat{f}(A)| \leq \Psi_j.$$

PROOF. The proof relies on a pair of VC inequalities (see [4], Chapter 3) that bound the *relative* deviation of true and empirical probabilities. For the collection $\mathcal{A}_j$ with cardinality $2^{jd}$, the relative VC inequalities imply that for any $\varepsilon > 0$, with probability $> 1 - 8 \cdot 2^{jd} e^{-n\varepsilon^2/4}$, $\forall A \in \mathcal{A}_j$ both

$$P(A) - \widehat{P}(A) \leq \varepsilon \sqrt{P(A)} \quad \text{and} \quad \widehat{P}(A) - P(A) \leq \varepsilon \sqrt{\widehat{P}(A)}.$$

Also, observe that

(4) $\qquad \widehat{P}(A) \leq P(A) + \varepsilon \sqrt{\widehat{P}(A)} \quad \Longrightarrow \quad \widehat{P}(A) \leq 2 \max(P(A), 2\varepsilon^2)$

and

(5) $\qquad P(A) \leq \widehat{P}(A) + \varepsilon \sqrt{P(A)} \quad \Longrightarrow \quad P(A) \leq 2 \max(\widehat{P}(A), 2\varepsilon^2).$

To understand statement (4), consider the following two cases: (i) If $\widehat{P}(A) \leq 4\varepsilon^2$, the statement is obvious; (ii) if $\widehat{P}(A) > 4\varepsilon^2$, this gives a bound on $\varepsilon$ which implies $\widehat{P}(A) \leq P(A) + \widehat{P}(A)/2 \Longrightarrow \widehat{P}(A) \leq 2P(A)$. Statement (5) follows similarly. Therefore, using (5) we get, with probability $> 1 - 8 \cdot 2^{jd} e^{-n\varepsilon^2/4}$, $\forall A \in \mathcal{A}_j$,

$$|P(A) - \widehat{P}(A)| \leq \varepsilon \sqrt{2 \max(\widehat{P}(A), 2\varepsilon^2)}.$$

Setting $\varepsilon = \sqrt{4 \log(2^{jd} 8/\delta_j)/n}$, $\delta_j = \delta 2^{-(j+1)}$ and applying union bound, we have with probability $> 1 - \delta$, for all $j \geq 0$ and all cells $A \in \mathcal{A}_j$

$$|P(A) - \widehat{P}(A)| \leq \sqrt{8 \frac{\log(2^{j(d+1)} 16/\delta)}{n} \max\left(\widehat{P}(A), 8 \frac{\log(2^{j(d+1)} 16/\delta)}{n}\right)}.$$

The result follows by dividing both sides by $\mu(A)$. □

The next lemma states how the density deviation bound or penalty $\Psi_j$ scales with resolution $j$ and number of observations $n$.

LEMMA A.2. *There exist constants $c_3, c_4 \equiv c_4(f_{\max}, d) > 0$ such that if $j \equiv j(n)$ satisfies $2^j = O((n/\log n)^{1/d})$, then for all $n$, with probability at least $1 - 1/n$,*

$$c_3 \sqrt{2^{jd} \frac{\log n}{n}} \leq \Psi_j \leq c_4 \sqrt{2^{jd} \frac{\log n}{n}}.$$

PROOF. We first derive the lower bound. Observe that since the total empirical probability mass is 1, we have

$$1 = \sum_{A \in \mathcal{A}_j} \widehat{P}(A) \leq \max_{A \in \mathcal{A}_j} \widehat{P}(A) \times |\mathcal{A}_j| = \max_{A \in \mathcal{A}_j} \frac{\widehat{P}(A)}{\mu(A)} = \max_{A \in \mathcal{A}_j} \widehat{f}(A).$$

Using this along with $\delta = 1/n$, $j \geq 0$ and $\mu(A) = 2^{-jd}$, we get

$$\Psi_j \geq \sqrt{2^{jd} 8 \frac{\log 16n}{n}}.$$

To get the upper bound, using statement (4) from the proof of Lemma A.1, we have, with probability $> 1 - 8 \cdot 2^{jd} e^{-n\varepsilon^2/4}$, for all $A \in \mathcal{A}_j$, $\widehat{P}(A) \leq 2\max(P(A), 2\varepsilon^2)$. Setting $\varepsilon = \sqrt{4\log(2^{jd}8/\delta_j)/n}$, $\delta_j = \delta 2^{-(j+1)}$ and applying union bound, we have, with probability $> 1 - \delta$, for all $j \geq 0$ and all $A \in \mathcal{A}_j$,

$$\widehat{P}(A) \leq 2\max\left(P(A), 8\frac{\log(2^{j(d+1)}16/\delta)}{n}\right).$$

Dividing by $\mu(A) = 2^{-jd}$ and using the density bound $f_{\max}$, we get a bound on $\max_{A \in \mathcal{A}_j} \widehat{f}(A)$, which implies that, with probability $> 1 - \delta$,

$$\Psi_j \leq \sqrt{2^{jd} 8 \frac{\log(2^{j(d+1)}16/\delta)}{n} \cdot 2\max\left(f_{\max}, 2^{jd}8\frac{\log(2^{j(d+1)}16/\delta)}{n}\right)}.$$

And using $\delta = 1/n$ and $2^j = O((n/\log n)^{1/d})$, we get

$$\Psi_j \leq c_4(f_{\max}, d)\sqrt{2^{jd} \frac{\log n}{n}}. \qquad \square$$

We now analyze the performance of the plug-in histogram-based level set estimator proposed in (1), and establish the following lemma that bounds its Hausdorff error. The first term denotes the estimation error while the second term that is proportional to the sidelength of a cell ($2^{-j}$) reflects the approximation error. We would like to point out that some arguments in the proofs hold for $s_n$ large enough. This implies that some of the constants in our proofs will depend on $\{s_i\}_{i=1}^\infty$, the exact form that the sequence $s_n$ takes (but not on $n$). However, we omit this dependence for simplicity.

LEMMA A.3. *Consider densities satisfying assumptions* [A1] *and* [B]. *If* $j \equiv j(n)$ *is such that* $2^j = O(s_n^{-1}(n/\log n)^{1/d})$, *where* $s_n$ *is a monotone diverging sequence, and* $n \geq n_0 \equiv n_0(f_{\max}, d, \delta_1, \varepsilon_o, C_1, \alpha)$, *then with probability at least* $1 - 3/n$,

$$d_\infty(\widehat{G}_j, G_\gamma^*) \leq \max(2C_3 + 3, 8\sqrt{d}\varepsilon_o^{-1})\left[\left(\frac{\Psi_j}{C_1}\right)^{1/\alpha} + \sqrt{d}2^{-j}\right].$$

PROOF. Let $J_0 = \lceil \log_2 4\sqrt{d}/\varepsilon_o \rceil$, where $\varepsilon_o$ is as defined in assumption [B]. Also, define

$$\varepsilon_j := \left[\left(\frac{\Psi_j}{C_1}\right)^{1/\alpha} + \sqrt{d}2^{-j}\right].$$

Consider the following two cases:

I. $j < J_0$. For this case, since the domain $\mathcal{X} = [0,1]^d$, we use the trivial bound
$$d_\infty(\widehat{G}_j, G_\gamma^*) \leq \sqrt{d} \leq 2^{J_0}(\sqrt{d}2^{-j}) \leq 8\sqrt{d}\varepsilon_o^{-1}\varepsilon_j.$$
The last step follows by choice of $J_0$ and since $\Psi_j, C_1 > 0$.

II. $j \geq J_0$. Observe that assumption [B] implies that $G_\gamma^*$ is not empty since $G_\gamma^* \supseteq \mathcal{I}_\varepsilon(G_\gamma^*) \neq \varnothing$ for $\varepsilon \leq \varepsilon_o$. We will show that for large enough $n$, with high probability, $\widehat{G}_j \cap G_\gamma^* \neq \varnothing$ for $j \geq J_0$, and hence $\widehat{G}_j$ is not empty. Thus, the Hausdorff error is given as

$$d_\infty(\widehat{G}_j, G_\gamma^*) = \max\Big\{\sup_{x \in G_\gamma^*} \rho(x, \widehat{G}_j), \sup_{x \in \widehat{G}_j} \rho(x, G_\gamma^*)\Big\}, \tag{6}$$

and we need bounds on the two terms in the right-hand side.

To prove that $\widehat{G}_j$ is not empty and obtain bounds on the two terms in the Hausdorff error, we establish a proposition and corollary. In the following analysis, if $G = \varnothing$, then we define $\sup_{x \in G} g(x) = 0$ for any function $g(\cdot)$. The proposition establishes that for large enough $n$, with high probability, all points whose distance to the boundary $\partial G_\gamma^*$ is greater than $\varepsilon_j$ are correctly excluded or included in the level set estimate.

PROPOSITION 2. *If $j \equiv j(n)$ is such that $2^j = O(s_n^{-1}(n/\log n)^{1/d})$, and $n \geq n_1(f_{\max}, d, \delta_1)$, then with probability at least $1 - 2/n$,*

$$\sup_{x \in \widehat{G}_j \Delta G_\gamma^*} \rho(x, \partial G_\gamma^*) \leq \left(\frac{\Psi_j}{C_1}\right)^{1/\alpha} + \sqrt{d}2^{-j} = \varepsilon_j.$$

PROOF. If $\widehat{G}_j \Delta G_\gamma^* = \varnothing$, then $\sup_{x \in \widehat{G}_j \Delta G_\gamma^*} \rho(x, \partial G_\gamma^*) = 0$ by definition, and the result of the proposition holds. If $\widehat{G}_j \Delta G_\gamma^* \neq \varnothing$, consider $x \in \widehat{G}_j \Delta G_\gamma^*$. Let $A_x \in \mathcal{A}_j$ denote the cell containing $x$ at resolution $j$. Consider the following two cases:

(i) $A_x \cap \partial G_\gamma^* \neq \varnothing$. This implies that $\rho(x, \partial G_\gamma^*) \leq \sqrt{d}2^{-j}$.

(ii) $A_x \cap \partial G_\gamma^* = \varnothing$. Since $x \in \widehat{G}_j \Delta G_\gamma^*$, it is erroneously included or excluded from the level set estimate $\widehat{G}_j$. Therefore, if $\bar{f}(A_x) \geq \gamma$, then $\widehat{f}(A_x) < \gamma$ and if $\bar{f}(A_x) < \gamma$, then $\widehat{f}(A_x) \geq \gamma$. This implies that $|\gamma - \bar{f}(A_x)| \leq |\bar{f}(A_x) - \widehat{f}(A_x)|$. Using Lemma A.1, we get $|\gamma - \bar{f}(A_x)| \leq \Psi_j$ with probability at least $1 - \delta$.

Now let $x_1$ be any point in $A_x$ such that $|\gamma - f(x_1)| \leq |\gamma - \bar{f}(A_x)|$. (Notice that at least one such point must exist in $A_x$ since this cell does not intersect the boundary.) As argued above, $|\gamma - \bar{f}(A_x)| \leq \Psi_j$ with probability at least $1 - 1/n$ (for $\delta = 1/n$). Using Lemma A.2, for resolutions satisfying $2^j = O(s_n^{-1}(n/\log n)^{1/d})$ and for large enough $n \geq n_1(f_{\max}, d, \delta_1)$, $\Psi_j \leq \delta_1$; hence, $|\gamma - f(x_1)| \leq \delta_1$, with probability at least $1 - 1/n$. Thus, the density regularity assumption [A1] holds at $x_1$ with probability $> 1 - 2/n$, and we have

$$\rho(x_1, \partial G_\gamma^*) \leq \left(\frac{|\gamma - f(x_1)|}{C_1}\right)^{1/\alpha} \leq \left(\frac{|\gamma - \bar{f}(A_x)|}{C_1}\right)^{1/\alpha} \leq \left(\frac{\Psi_j}{C_1}\right)^{1/\alpha}.$$

Since $x, x_1 \in A_x$,

$$\rho(x, \partial G_\gamma^*) \leq \rho(x_1, \partial G_\gamma^*) + \sqrt{d}2^{-j} \leq \left(\frac{\Psi_j}{C_1}\right)^{1/\alpha} + \sqrt{d}2^{-j}.$$

So for both cases, if $j \equiv j(n)$ is such that $2^j = O(s_n^{-1}(n/\log n)^{1/d})$, and $n \geq n_1(f_{\max}, d, \delta_1)$, then with probability at least $1 - 2/n$, $\forall x \in \widehat{G}_j \Delta G_\gamma^*$, $\rho(x, \partial G_\gamma^*) \leq (\Psi_j/C_1)^{1/\alpha} + \sqrt{d}2^{-j} = \varepsilon_j$. □

Based on Proposition 2, the following corollary argues that for large enough $n$ and $j \geq J_0 = \lceil \log_2 4\sqrt{d}/\varepsilon_o \rceil$, with high probability, all points within the inner cover $\mathcal{I}_{2\varepsilon_j}(G_\gamma^*)$ that are at a distance greater than $\varepsilon_j$ are correctly included in the level set estimate; hence, they lie in $\widehat{G}_j \cap G_\gamma^*$. This also implies that $\widehat{G}_j$ is not empty.

COROLLARY 1. *Recall assumption* [B] *and denote the inner cover of $G_\gamma^*$ with $2\varepsilon_j$-balls, $\mathcal{I}_{2\varepsilon_j}(G_\gamma^*) \equiv \mathcal{I}_{2\varepsilon_j}$ for simplicity. If $j \equiv j(n)$ is such that $2^j = O(s_n^{-1}(n/\log n)^{1/d})$, $j \geq J_0$, and $n \geq n_0 \equiv n_0(f_{\max}, d, \delta_1, \varepsilon_o, C_1, \alpha)$, then with probability at least $1 - 3/n$,*

$$\widehat{G}_j \neq \varnothing \quad \text{and} \quad \sup_{x \in \mathcal{I}_{2\varepsilon_j}} \rho(x, \widehat{G}_j \cap G_\gamma^*) \leq \varepsilon_j.$$

PROOF. Observe that for $j \geq J_0$, $2\sqrt{d}2^{-j} \leq 2\sqrt{d}2^{-J_0} \leq \varepsilon_o/2$. By Lemma A.2, for resolutions satisfying $2^j = O(s_n^{-1}(n/\log n)^{1/d})$, and for large enough $n \geq n_2(\varepsilon_o, f_{\max}, C_1, \alpha)$, $2(\Psi_j/C_1)^{1/\alpha} \leq \varepsilon_o/2$, with probability at least $1 - 1/n$. Therefore, for resolutions satisfying $2^j = O(s_n^{-1}(n/\log n)^{1/d})$ and $j \geq J_0$, and for $n \geq n_2$, with probability at least $1 - 1/n$, $2\varepsilon_j \leq \varepsilon_o$ and hence $\mathcal{I}_{2\varepsilon_j} \neq \varnothing$.

Now consider any $2\varepsilon_j$-ball in $\mathcal{I}_{2\varepsilon_j}$. Then the distance of all points in the interior of the concentric $\varepsilon_j$-ball from the boundary of $\mathcal{I}_{2\varepsilon_j}$, and hence from the boundary of $G_\gamma^*$, is greater than $\varepsilon_j$. As per Proposition 2, for $n \geq n_0 = \max(n_1, n_2)$ with probability $> 1 - 3/n$, none of these points can lie in $\widehat{G}_j \Delta G_\gamma^*$; hence, they must lie in $\widehat{G}_j \cap G_\gamma^*$ since they are in $\mathcal{I}_{2\varepsilon_j} \subseteq G_\gamma^*$. Thus, $\widehat{G}_j \neq \varnothing$, and for all $x \in \mathcal{I}_{2\varepsilon_j}$, $\rho(x, \widehat{G}_j \cap G_\gamma^*) \leq \varepsilon_j$. □

We now resume the proof of Lemma A.3, case II. Assume the conclusions of Proposition 2 and Corollary 1 hold. Thus, all the following statements hold for resolutions satisfying $2^j = O(s_n^{-1}(n/\log n)^{1/d})$, $j \geq J_0$ and $n \geq n_0 \equiv n_0(f_{\max}, d, \delta_1, \varepsilon_o, C_1, \alpha)$, with probability at least $1 - 3/n$. Since $G_\gamma^*$ and $\widehat{G}_j$ are nonempty sets, we now bound the two terms that contribute to the Hausdorff error. To bound the term $\sup_{x \in \widehat{G}_j} \rho(x, G_\gamma^*)$, observe that

$$\sup_{x \in \widehat{G}_j} \rho(x, G_\gamma^*) = \sup_{x \in \widehat{G}_j \setminus G_\gamma^*} \rho(x, G_\gamma^*) = \sup_{x \in \widehat{G}_j \setminus G_\gamma^*} \rho(x, \partial G_\gamma^*)$$
(7)
$$\leq \sup_{x \in \widehat{G}_j \Delta G_\gamma^*} \rho(x, \partial G_\gamma^*) \leq \varepsilon_j,$$

where the last step follows from Proposition 2.

To bound the term $\sup_{x \in G_\gamma^*} \rho(x, \widehat{G}_j)$, we recall assumption [B] which states that the boundary points of $G_\gamma^*$ are $O(\varepsilon_j)$ from the inner cover $\mathcal{I}_{2\varepsilon_j}(G_\gamma^*)$, and we use Corollary 1 to bound the distance of the inner cover from $\widehat{G}_j$ as follows:

$$\sup_{x \in G_\gamma^*} \rho(x, \widehat{G}_j) \leq \sup_{x \in G_\gamma^*} \rho(x, \widehat{G}_j \cap G_\gamma^*)$$
(8)
$$= \max\left\{\sup_{x \in \mathcal{I}_{2\varepsilon_j}} \rho(x, \widehat{G}_j \cap G_\gamma^*), \sup_{x \in G_\gamma^* \setminus \mathcal{I}_{2\varepsilon_j}} \rho(x, \widehat{G}_j \cap G_\gamma^*)\right\}$$
$$\leq \max\left\{\varepsilon_j, \sup_{x \in G_\gamma^* \setminus \mathcal{I}_{2\varepsilon_j}} \rho(x, \widehat{G}_j \cap G_\gamma^*)\right\},$$

where the last step follows from Corollary 1.

Now consider any $x \in G_\gamma^* \setminus \mathcal{I}_{2\varepsilon_j}$. By the triangle inequality, $\forall y \in \partial G_\gamma^*$ and $\forall z \in \mathcal{I}_{2\varepsilon_j}$,

$$\rho(x, \widehat{G}_j \cap G_\gamma^*) \leq \rho(x, y) + \rho(y, z) + \rho(z, \widehat{G}_j \cap G_\gamma^*)$$
$$\leq \rho(x, y) + \rho(y, z) + \sup_{z' \in \mathcal{I}_{2\varepsilon_j}} \rho(z', \widehat{G}_j \cap G_\gamma^*)$$
$$\leq \rho(x, y) + \rho(y, z) + \varepsilon_j,$$

where the last step follows from Corollary 1. This implies that, $\forall y \in \partial G_\gamma^*$,

$$\rho(x, \widehat{G}_j \cap G_\gamma^*) \leq \rho(x, y) + \inf_{z \in \mathcal{I}_{2\varepsilon_j}} \rho(y, z) + \varepsilon_j$$
$$= \rho(x, y) + \rho(y, \mathcal{I}_{2\varepsilon_j}) + \varepsilon_j$$

$$\leq \rho(x,y) + \sup_{y' \in \partial G_\gamma^*} \rho(y', \mathcal{I}_{2\varepsilon_j}) + \varepsilon_j$$

$$\leq \rho(x,y) + 2C_3\varepsilon_j + \varepsilon_j,$$

where the last step invokes assumption [B]. This, in turn, implies

$$\rho(x, \widehat{G}_j \cap G_\gamma^*) \leq \inf_{y \in \partial G_\gamma^*} \rho(x,y) + (2C_3+1)\varepsilon_j \leq 2\varepsilon_j + (2C_3+1)\varepsilon_j.$$

The second step is true for $x \in G_\gamma^* \setminus \mathcal{I}_{2\varepsilon_j}$, because if it was not true then $\forall y \in \partial G_\gamma^*$, $\rho(x,y) > 2\varepsilon_j$; hence, there exists a closed $2\varepsilon_j$-ball around $x$ that is in $G_\gamma^*$. This contradicts the fact that $x \notin \mathcal{I}_{2\varepsilon_j}$. Therefore, we have

$$\sup_{x \in G_\gamma^* \setminus \mathcal{I}_{2\varepsilon_j}} \rho(x, \widehat{G}_j \cap G_\gamma^*) \leq (2C_3+3)\varepsilon_j.$$

And going back to (8), we get

(9) $$\sup_{x \in G_\gamma^*} \rho(x, \widehat{G}_j) \leq (2C_3+3)\varepsilon_j.$$

From (7) and (9), we have that for all densities satisfying assumptions [A1] and [B], if $j \equiv j(n)$ is such that $2^j = O(s_n^{-1}(n/\log n)^{1/d})$, $j \geq J_0$ and $n \geq n_0 \equiv n_0(f_{\max}, d, \delta_1, \varepsilon_o, C_1, \alpha)$, then with probability $> 1 - 3/n$,

$$d_\infty(\widehat{G}_j, G_\gamma^*) = \max\Big\{\sup_{x \in G_\gamma^*} \rho(x, \widehat{G}_j), \sup_{x \in \widehat{G}_j} \rho(x, G_\gamma^*)\Big\} \leq (2C_3+3)\varepsilon_j.$$

And addressing both case I ($j < J_0$) and case II ($j \geq J_0$), we finally have that for all densities satisfying assumptions [A1] and [B], if $j \equiv j(n)$ is such that $2^j = O(s_n^{-1}(n/\log n)^{1/d})$, and $n \geq n_0 \equiv n_0(f_{\max}, d, \delta_1, \varepsilon_o, C_1, \alpha)$, then with probability $> 1 - 3/n$,

$$d_\infty(\widehat{G}_j, G_\gamma^*) \leq \max(2C_3+3, 8\sqrt{d}\varepsilon_o^{-1})\varepsilon_j.$$

This completes the proof of Lemma A.3. □

We now establish the result of Theorem 1. Since the local density regularity parameter $\alpha$ is known, the appropriate histogram resolution can be chosen as $2^{-j} \asymp s_n(n/\log n)^{-1/(d+2\alpha)}$. Let $\Omega$ denote the event such that the bounds of Lemma A.2 (with $\delta = 1/n$) and Lemma A.3 hold. Then for $n \geq n_0$, $P(\bar{\Omega}) \leq 4/n$, where $\bar{\Omega}$ denotes the complement of $\Omega$. For $n < n_0$, we can use the trivial inequality $P(\bar{\Omega}) \leq 1$. So we have, for all $n$, $P(\bar{\Omega}) \leq \max(4, n_0)\frac{1}{n} =: C'\frac{1}{n}$. Here $C' \equiv C'(f_{\max}, d, \delta_1, \varepsilon_o, C_1, \alpha)$. So $\forall f \in \mathcal{F}_1^*(\alpha)$, we have the following. (Explanation for each step is provided after the equations.)

$$\mathbb{E}[d_\infty(\widehat{G}_j, G_\gamma^*)] = P(\Omega)\mathbb{E}[d_\infty(\widehat{G}_j, G_\gamma^*)|\Omega] + P(\bar{\Omega})\mathbb{E}[d_\infty(\widehat{G}_j, G_\gamma^*)|\bar{\Omega}]$$

$$\leq \mathbb{E}[d_\infty(\widehat{G}_j, G_\gamma^*)|\Omega] + P(\bar{\Omega})\sqrt{d}$$

$$\leq \max(2C_3 + 3, 8\sqrt{d}\varepsilon_o^{-1}) \left[ \left(\frac{\Psi_j}{C_1}\right)^{1/\alpha} + \sqrt{d}2^{-j} \right] + C'\frac{\sqrt{d}}{n}$$

$$\leq C \max\left\{ \left(2^{jd}\frac{\log n}{n}\right)^{1/(2\alpha)}, 2^{-j}, \frac{1}{n} \right\}$$

$$\leq C \max\left\{ s_n^{-d/2\alpha}\left(\frac{n}{\log n}\right)^{-1/(d+2\alpha)}, s_n\left(\frac{n}{\log n}\right)^{-1/(d+2\alpha)}, \frac{1}{n} \right\}$$

$$\leq C s_n \left(\frac{n}{\log n}\right)^{-1/(d+2\alpha)}.$$

Here $C \equiv C(C_1, C_3, \varepsilon_o, f_{\max}, \delta_1, d, \alpha)$. The second step follows by observing the trivial bounds $P(\Omega) \leq 1$ and $\mathbb{E}[d_\infty(\widehat{G}_j, G_\gamma^*)|\bar{\Omega}] \leq \sqrt{d}$ since the domain $\mathcal{X} = [0,1]^d$. The third step follows from Lemma A.3 and the fourth one using Lemma A.2. The fifth step follows since the chosen resolution $2^{-j} \asymp s_n(n/\log n)^{-1/(d+2\alpha)}$.

## APPENDIX B: PROOF OF THEOREM 2

To analyze the resolution chosen by the complexity penalized procedure of (3) based on the vernier, we first establish two results regarding the vernier. Using Lemma A.1, we have the following corollary that bounds the deviation of true and empirical vernier.

COROLLARY B.1. *Consider $0 < \delta < 1$. With probability at least $1 - \delta$, the following is true for all $j \geq 0$:*

$$|\mathcal{V}_{\gamma,j} - \widehat{\mathcal{V}}_{\gamma,j}| \leq \Psi_{j'}.$$

PROOF. Let $A_0 \in \mathcal{A}_j$ denote the cell achieving the minimum defining $\mathcal{V}_{\gamma,j}$ and $A_1 \in \mathcal{A}_j$ denote the cell achieving the minimum defining $\widehat{\mathcal{V}}_{\gamma,j}$. Also, let $A'_{00}$ and $A'_{10}$ denote the subcells at resolution $j'$ within $A_0$ and $A_1$, respectively, that have maximum average density deviation from $\gamma$. Similarly, let $A'_{01}$ and $A'_{11}$ denote the subcells at resolution $j'$ within $A_0$ and $A_1$, respectively, that have maximum empirical density deviation from $\gamma$. Then, we have

$$\begin{aligned}
\mathcal{V}_{\gamma,j} - \widehat{\mathcal{V}}_{\gamma,j} &= |\gamma - \bar{f}(A'_{00})| - |\gamma - \widehat{f}(A'_{11})| \\
&\leq |\gamma - \bar{f}(A'_{10})| - |\gamma - \widehat{f}(A'_{11})| \leq |\bar{f}(A'_{10}) - \widehat{f}(A'_{11})| \\
&\leq \max\{\bar{f}(A'_{10}) - \widehat{f}(A'_{10}), \widehat{f}(A'_{11}) - \bar{f}(A'_{11})\} \\
&\leq \max_{A \in \mathcal{A}_{j'}} |\bar{f}(A) - \widehat{f}(A)| \leq \Psi_{j'}.
\end{aligned}$$

The first inequality invokes definition of $A_0$, the third inequality invokes definitions of the subcells $A'_{10}$, $A'_{11}$ and the last one follows from Lemma A.1. The bound on $\widehat{\mathcal{V}}_{\gamma,j} - \mathcal{V}_{\gamma,j}$ follows similarly. $\square$

The second result establishes that the vernier is sensitive to the resolution and density regularity.

LEMMA B.1. *Consider densities satisfying assumptions* [A] *and* [B]. *Recall that* $j' = \lfloor j + \log_2 s_n \rfloor$, *where* $s_n$ *is a monotone diverging sequence. There exists* $C \equiv C(C_2, f_{\max}, \delta_2, \alpha) > 0$ *such that if* $n$ *is large enough so that* $s_n > 8\max(3\varepsilon_o^{-1}, 28, 12C_3)\sqrt{d}$, *then for all* $j \geq 0$,

$$\min(\delta_1, C_1)2^{-j'\alpha} \leq \mathcal{V}_{\gamma,j} \leq C(\sqrt{d}2^{-j})^\alpha.$$

PROOF. We first establish the upper bound. Recall assumption [A] and consider the cell $A_0 \in \mathcal{A}_j$ that contains the point $x_0$. Then, $A_0 \cap \partial G^*_\gamma \neq \varnothing$. Let $A'_0$ denote the subcell at resolution $j'$ within $A_0$ that has maximum average density deviation from $\gamma$. Consider the following two cases:

(i) If the resolution is high enough so that $\sqrt{d}2^{-j} \leq \delta_2$, then the density regularity assumption [A2] holds $\forall x \in A_0$ since $A_0 \subset B(x_0, \delta_2)$, the $\delta_2$-ball around $x_0$. The same holds also for the subcell $A'_0$. Hence,

$$|\gamma - \bar{f}(A'_0)| \leq C_2(\sqrt{d}2^{-j})^\alpha.$$

(ii) If the resolution is not high enough and $\sqrt{d}2^{-j} > \delta_2$, use the following trivial bound: $|\gamma - \bar{f}(A'_0)| \leq f_{\max} \leq \frac{f_{\max}}{\delta_2^\alpha}(\sqrt{d}2^{-j})^\alpha$.

Hence, we can say for all $j$ there exists a cell $A_0 \in \mathcal{A}_j$ such that

$$\max_{A' \in \mathcal{A}_{j'} \cap A_0} |\gamma - \bar{f}(A')| = |\gamma - \bar{f}(A'_0)| \leq \max\left(C_2, \frac{f_{\max}}{\delta_2^\alpha}\right)(\sqrt{d}2^{-j})^\alpha.$$

This yields the upper bound on the vernier, $\mathcal{V}_{\gamma,j} \leq C(\sqrt{d}2^{-j})^\alpha$, where $C \equiv C(C_2, f_{\max}, \delta_2, \alpha)$.

For the lower bound, consider any cell $A \in \mathcal{A}_j$. We will show that the level set regularity assumption [B] implies that for large enough $n$ (so that the sidelength $2^{-j'}$ is small enough), the boundary does not intersect all subcells at resolution $j'$ within the cell $A$ at resolution $j$. In fact, there exists at least one subcell $A'_1 \in A \cap \mathcal{A}_{j'}$ such that $\forall x \in A'_1$,

$$\rho(x, \partial G^*_\gamma) \geq 2^{-j'}.$$

We establish this statement formally later on, but for now assume that it holds. The local density regularity condition [A] now gives that for all $x \in A'_1$, $|\gamma - f(x)| \geq \min(\delta_1, C_1 2^{-j'\alpha}) \geq \min(\delta_1, C_1)2^{-j'\alpha}$. So we have

$$\max_{A' \in A \cap \mathcal{A}_{j'}} |\gamma - \bar{f}(A')| \geq |\gamma - \bar{f}(A'_1)| \geq \min(\delta_1, C_1)2^{-j'\alpha}.$$

Since this is true for any $A \in \mathcal{A}_j$, in particular, this is true for the cell achieving the minimum defining $\mathcal{V}_{\gamma,j}$. Hence, the lower bound on the vernier $\mathcal{V}_{\gamma,j}$ follows.

We now formally prove that the level set regularity assumption [B] implies that for large enough $n$ (so that $s_n > 8\max(3\varepsilon_o^{-1}, 28, 12C_3)\sqrt{d}$), $\exists A_1' \in A \cap \mathcal{A}_{j'}$ such that $\forall x \in A_1'$,

$$\rho(x, \partial G_\gamma^*) \geq 2^{-j'}.$$

Observe that if we consider any cell at resolution $j'' := j' - 2$ that does not intersect the boundary $\partial G_\gamma^*$, then it contains a cell at resolution $j'$ that is greater than $2^{-j'}$ away from the boundary. Thus, it suffices to show that for large enough $n$ [so that $s_n > 8\max(3\varepsilon_o^{-1}, 28, 12C_3)\sqrt{d}$], $\exists A'' \in A \cap \mathcal{A}_{j''}$ such that $A'' \cap \partial G_\gamma^* = \varnothing$. We prove the last statement by contradiction. Suppose that for $s_n > 8\max(3\varepsilon_o^{-1}, 28, 12C_3)\sqrt{d}$, all subcells in $A$ at resolution $j''$ intersect the boundary $\partial G_\gamma^*$. Let $\varepsilon = 3\sqrt{d}2^{-j''}$. Then,

$$\varepsilon = 3\sqrt{d}2^{-j''} = 12\sqrt{d}2^{-j'} < \frac{24\sqrt{d}}{s_n}2^{-j} \leq \frac{24\sqrt{d}}{s_n} \leq \varepsilon_o,$$

where the last step follows since $s_n \geq 24\sqrt{d}\varepsilon_o^{-1}$. By choice of $\varepsilon$, every closed $\varepsilon$-ball in $A$ must contain an entire subcell at resolution $j''$ and in fact must contain an open neighborhood around that subcell. Since the boundary intersects all subcells at resolution $j''$, this implies that every closed $\varepsilon$-ball in $A$ contains a boundary point and in fact contains an open neighborhood around that boundary point. Thus, (i) every closed $\varepsilon$-ball in $A$ contains points not in $G_\gamma^*$, and hence cannot lie in $\mathcal{I}_\varepsilon(G_\gamma^*)$. Also, observe that since all subcells in $A$ at resolution $j''$ intersect the boundary of $G_\gamma^*$, (ii) there exists a boundary point $x_1$ that is within $\sqrt{d}2^{-j''}$ of the center of cell $A$. From (i) and (ii) it follows that

$$\rho(x_1, \mathcal{I}_\varepsilon(G_\gamma^*)) \geq \frac{2^{-j}}{2} - \sqrt{d}2^{-j''} - 2\varepsilon = \frac{2^{-j}}{2} - 28\sqrt{d}2^{-j'}$$

$$> 2^{-j}\left(\frac{1}{2} - \frac{56\sqrt{d}}{s_n}\right) > \frac{2^{-j}}{4},$$

where the last step follows since $s_n > 224\sqrt{d}$. However, assumption [B] implies that for $\varepsilon \leq \varepsilon_o$,

$$\rho(x_1, \mathcal{I}_\varepsilon(G_\gamma^*)) \leq C_3\varepsilon = 3C_3\sqrt{d}2^{-j''} = 12C_3\sqrt{d}2^{-j'} \leq \frac{24C_3\sqrt{d}2^{-j}}{s_n} \leq \frac{2^{-j}}{4},$$

where the last step follows since $s_n > 96C_3\sqrt{d}$, and we have a contradiction.

This completes the proof of Lemma B.1. $\square$

We are now ready to prove Theorem 2. To analyze the resolution $\widehat{j}$ chosen by (3), we first derive upper bounds on $\mathcal{V}_{\gamma,\widehat{j}}$ and $\Psi_{\widehat{j}'}$ that effectively characterize the approximation error and estimation error, respectively. Thus, a bound on the vernier $\mathcal{V}_{\gamma,\widehat{j}}$ will imply that the chosen resolution $\widehat{j}$ cannot be too coarse, and a bound on the penalty will imply that the chosen resolution is not too fine. Using Corollary B.1 and (3), we have the following oracle inequality that holds with probability at least $1-\delta$:

$$\mathcal{V}_{\gamma,\widehat{j}} \leq \widehat{\mathcal{V}}_{\gamma,\widehat{j}} + \Psi_{\widehat{j}'} = \min_{0 \leq j \leq J}\{\widehat{\mathcal{V}}_{\gamma,j} + \Psi_{j'}\} \leq \min_{0 \leq j \leq J}\{\mathcal{V}_{\gamma,j} + 2\Psi_{j'}\}.$$

Lemma B.1 provides an upper bound on the vernier $\mathcal{V}_{\gamma,j}$, and Lemma A.2 provides an upper bound on the penalty $\Psi_{j'}$. We plug these bounds into the oracle inequality. Here $C$ may denote a different constant from line to line. With probability at least $1 - 2/n$ (with $\delta = 1/n$),

$$\mathcal{V}_{\gamma,\widehat{j}} \leq \widehat{\mathcal{V}}_{\gamma,\widehat{j}} + \Psi_{\widehat{j}'} \leq C \min_{0 \leq j \leq J}\left\{\max\left(2^{-j\alpha}, \sqrt{2^{jd}s_n^d \frac{\log n}{n}}\right)\right\}$$

$$\leq C s_n^{d\alpha/(d+2\alpha)}\left(\frac{n}{\log n}\right)^{-\alpha/(d+2\alpha)}.$$

Here $C \equiv C(C_2, f_{\max}, \delta_2, d, \alpha)$. The first step uses the definition of $j'$, and the second step follows by balancing the two terms for optimal resolution $j^*$ given by $2^{-j^*} \asymp s_n^{d/(d+2\alpha)}(n/\log n)^{-1/(d+2\alpha)}$. This establishes the desired bounds on $\mathcal{V}_{\gamma,\widehat{j}}$ and $\Psi_{\widehat{j}'}$.

Now, using Lemma B.1 and the definition of $j'$, we have the following upper bound on the chosen sidelength. For $s_n > 8\max(3\varepsilon_o^{-1}, 28, 12C_3)\sqrt{d}$,

$$2^{-\widehat{j}} \leq s_n 2^{-\widehat{j}'} \leq s_n\left(\frac{\mathcal{V}_{\gamma,\widehat{j}}}{\min(\delta_1, C_1)}\right)^{1/\alpha} \leq c_2 s_n s_n^{d/(d+2\alpha)}\left(\frac{n}{\log n}\right)^{-1/(d+2\alpha)},$$

where $c_2 \equiv c_2(C_1, C_2, f_{\max}, \delta_1, \delta_2, d, \alpha) > 0$. Also notice that since $2^J \asymp s_n^{-1}(n/\log n)^{1/d}$, we have $2^{j'} \leq s_n 2^j \leq s_n 2^J \asymp (n/\log n)^{1/d}$, and thus $j'$ satisfies the condition of Lemma A.2. Therefore, using Lemma A.2, we get a lower bound on the sidelength. With probability at least $1 - 2/n$,

$$2^{-\widehat{j}} > \frac{s_n}{2} 2^{-\widehat{j}'} \geq \frac{s_n}{2}\left(\frac{\Psi_{\widehat{j}'}^2}{c_3^2}\frac{n}{\log n}\right)^{-1/d} \geq c_1 s_n^{d/(d+2\alpha)}\left(\frac{n}{\log n}\right)^{-1/(d+2\alpha)},$$

where $c_1 \equiv c_1(C_2, f_{\max}, \delta_2, d, \alpha) > 0$. So we have for $s_n > 8\max(3\varepsilon_o^{-1}, 28, 12C_3)\sqrt{d}$, with probability at least $1 - 2/n$,

$$(10) \quad c_1 s_n^{d/(d+2\alpha)}\left(\frac{n}{\log n}\right)^{-1/(d+2\alpha)} \leq 2^{-\widehat{j}} \leq c_2 s_n s_n^{d/(d+2\alpha)}\left(\frac{n}{\log n}\right)^{-1/(d+2\alpha)},$$

where $c_1 \equiv c_1(C_2, f_{\max}, \delta_2, d, \alpha) > 0$ and $c_2 \equiv c_2(C_1, C_2, f_{\max}, \delta_1, \delta_2, d, \alpha) > 0$. Hence, the automatically chosen resolution behaves as desired.

Now we can invoke Lemma A.3 to derive the rate of convergence for the Hausdorff error. Consider large enough $n \geq n_1(C_3, \varepsilon_o, d)$ so that $s_n > 8\max(3\varepsilon_o^{-1}, 28, 12C_3)\sqrt{d}$. Also, recall that the condition of Lemma A.3 requires that $n \geq n_0(f_{\max}, d, \delta_1, \varepsilon_o, C_1, \alpha)$. Pick $n \geq \max(n_0, n_1)$ and let $\Omega$ denote the event such that the bound of Lemma A.3 and the upper and lower bounds on the chosen resolution in (10) hold. Then, we have $P(\bar{\Omega}) \leq 5/n$. For $n < \max(n_0, n_1)$, we can use the trivial inequality $P(\bar{\Omega}) \leq 1$. So we have, for all $n$, $P(\bar{\Omega}) \leq \max(5, \max(n_0, n_1))\frac{1}{n} =: C\frac{1}{n}$. Here $C \equiv C(C_1, C_3, \varepsilon_o, f_{\max}, \delta_1, d, \alpha)$. So $\forall f \in \mathcal{F}_2^*(\alpha)$, we have the following. (Here $C$ may denote a different constant from line to line. Explanation for each step is provided after the equations.)

$$\mathbb{E}[d_\infty(\widehat{G}, G_\gamma^*)] = P(\Omega)\mathbb{E}[d_\infty(\widehat{G}, G_\gamma^*)|\Omega] + P(\bar{\Omega})\mathbb{E}[d_\infty(\widehat{G}, G_\gamma^*)|\bar{\Omega}]$$

$$\leq \mathbb{E}[d_\infty(\widehat{G}, G_\gamma^*)|\Omega] + P(\bar{\Omega})\sqrt{d}$$

$$\leq C\left[\left(\frac{\Psi_{\widehat{j}}}{C_1}\right)^{1/\alpha} + \sqrt{d}2^{-\widehat{j}} + \frac{\sqrt{d}}{n}\right]$$

$$\leq C\max\left\{\left(2^{\widehat{j}d}\frac{\log n}{n}\right)^{1/(2\alpha)}, 2^{-\widehat{j}}, \frac{1}{n}\right\}$$

$$\leq C\max\left\{s_n^{(-d^2/2\alpha)/d+2\alpha}\left(\frac{n}{\log n}\right)^{-1/(d+2\alpha)},\right.$$

$$\left. s_n s_n^{d/(d+2\alpha)}\left(\frac{n}{\log n}\right)^{-1/(d+2\alpha)}, \frac{1}{n}\right\}$$

$$\leq Cs_n^2\left(\frac{n}{\log n}\right)^{-1/(d+2\alpha)}.$$

Here $C \equiv C(C_1, C_2, C_3, \varepsilon_o, f_{\max}, \delta_1, \delta_2, d, \alpha)$. The second step follows by observing the trivial bounds $P(\Omega) \leq 1$ and since the domain $\mathcal{X} = [0,1]^d$, $\mathbb{E}[d_\infty(\widehat{G}, G_\gamma^*)|\bar{\Omega}] \leq \sqrt{d}$. The third step follows from Lemma A.3 and the fourth one from Lemma A.2. The fifth step follows using the upper and lower bounds established on $2^{-\widehat{j}}$ in (10).

**Acknowledgments.** The authors would like to thank Rui Castro and the anonymous reviewers for carefully reviewing the manuscript and providing valuable feedback that led to a much better presentation of the paper.

ADAPTIVE HAUSDORFF ESTIMATION OF DENSITY LEVEL SETS 23

A. SINGH  
R. NOWAK  
ELECTRICAL AND COMPUTER ENGINEERING  
UNIVERSITY OF WISCONSIN–MADISON  
3610 EH, 1415 ENGINEERING DRIVE  
MADISON, WISCONSIN 53726  
USA  
E-MAIL: singh@cae.wisc.edu

C. SCOTT  
ELECTRICAL ENGINEERING AND COMPUTER SCIENCE  
UNIVERSITY OF MICHIGAN–ANN ARBOR  
4229 EECS, 1301 BEAL AVENUE  
ANN ARBOR, MICHIGAN 48109  
USA